\documentclass[12pt]{amsart}

\usepackage{latexsym}
\usepackage{amsmath}
\usepackage{amssymb}
\usepackage{amsfonts}
\usepackage{a4}
\usepackage{alltt}
\usepackage{epsfig}

\def\cal{\mathcal}

\newcommand{\N}{{\mathbf N}}

\newtheorem{thm}{Theorem}[section] 
\newtheorem{prop}[thm]{Proposition} 
\newtheorem{lemma}[thm]{Lemma} 
 
\newtheorem{corr}[thm]{Corollary} 
 
\newtheorem{ex}[thm]{Example} 

\newtheorem{defn}[thm]{Definition}
\numberwithin{equation}{section}

\newcommand{\beq}{\begin{eqnarray}}
\newcommand{\eeq}{\end{eqnarray}}

\newcommand{\scbg}{semicentral bigroupoid}
\newcommand{\scbgs}{semicentral bigroupoids}

\newcommand{\rs}{rectangular structure}
\newcommand{\rss}{rectangular structures}

\title[Orderly Central Groupoids]
{Orderly Algorithm to enumerate central groupoids and their graphs}
\author[Tim Boykett]{Tim Boykett {\tt{tim@timesup.org}}
  Time's Up Research Department  and 
Department of Mathematics, Johannes--Kepler University, Linz, Austria}
\date{October 2003}

\begin{document}

\begin{abstract}

A graph has the unique path property $UPP_n$ if there is a
unique path of length $n$ between any ordered pair of nodes.

This paper reiterates Royle and MacKay's technique for constructing
orderly algorithms. We wish to use this technique to
enumerate all $UPP_2$ graphs of small orders $3^2$ and $4^2$.
We attempt to use the direct graph formalism and find that the
algorithm is inefficient.

We introduce a generalised problem and derive algebraic and
combinatoric structures with appropriate structure. We are
able to then design an orderly algorithm to determine all
$UPP_2$ graphs of order $3^2$, which runs fast enough. 
We hope to be able  to determine the $UPP_2$ graphs of order $4^2$
in the near future.
\end{abstract}

\maketitle


2000 Maths Subject: 05C30, 05C38, 05C50, 05C12

\section{Introduction}

Orderly algorithms are a concept for enumerating combinatoric objects:
trees, graphs, flocks and many other objects have been approached in
this way.
The essential hope is to efficiently enumerate the objects without generating
isomorphs: by cutting down the production, there is no need to filter the
isomorphs out afterwards.

Central groupoids were introduced by Evans in the late 1960s. Knuth later
showed that central groupoids are equivalent to a matrix problem proposed
by Hoffman: which 0-1 matrices $A$ exist such that $A^2 = J$ where $J$ is
the matrix consisting of all 1s. This is equivalent to the problem
of determining those directed graphs with exactly one path of length two
between any pair of nodes. Shader has shown that matrices of all possible
ranks exist. A series of papers by various researchers has investigated
generalisations of this problem and special solutions, but 
until recently 
exhaustive lists of examples of nontrivial order have not been published.
Recently a report listing all examples of order 9 has appeared. In the
following we will confirm these results, using techniques that
should be extendable to further listings.

This report is structured as follows. We introduce the orderly algorithm
as formulated by Royle, based upon the McKay formulation. This technique
is particularly suitable to structures that can be formulated in some
way as graphs, which will be our case.
We will then look at the central groupoids; some general theory and 
some existence results. Although there would seem to be a good connection between
the graph model of central groupoids and the orderly algorithm model
of Royle, this turns out to be inappropriate. We introduce a generalisation
of central groupoids and a derived combinatorial structure. It turns out that
this can be used in an orderly algorithm to efficiently enumerate examples.
We then use the theory connecting central groupoids to these generalisations
in order to filter out the desired examples.

The implementation in GAP \cite{GAP4} using the GRAPE \cite{soicher93}
package is discussed. GAP is a programming environment for group theory and other
algebraic programming. GRAPE is used for graph manipulation and 
to access the {\tt nauty} \cite{mckay90} package for graph automorphism.

We restrict ourselves to finite examples throughout.
The results are based upon related work with a distinct area
of application \cite{boykett04}.

\section{Orderly Algorithms}

A standard way to generate examples of combinatorial structures is
to generate examples in a way that can be proven to be exhaustive, and then
remove isomorphs from the list. This can be simple and effective,
but if we generate too many examples, then the effort required to filter
out the isomorphs can be overwhelming, owing to the quadratic
growth of the number of comparisons, a process that is often complex.

An alternative approach is to ensure that the
branching generation process does not generate isomorphic examples, thus
removing the filtering step of the algorithm described above.
With careful bookkeeping, this can be done.
Such generation processes have been called ``orderly'' in \cite{read78}.
McKay has developed a general structure for such
algorithms \cite{mckay98} and Royle has developed a simplified
algorithm \cite{royle98}, upon which we base ours.

Let $V$ be some set, $G$ a group acting upon $V$. We write $v^g$ for the
action of $g \in G$ on $v \in V$ and extend naturally
to actions on sets. We want to
find all subsets $X \subset V$ such that $P(X)$ is true
for some hereditary property $P$, but we want only
one example from each isomorphism class, with
isomorphism defined by $G$.

We require a function $\Theta$ such that
\beq
&& \Theta : 2^V \rightarrow 2^V \\
&&\Theta(X) \mbox{ is an orbit of } G_X \mbox{ on } X \\
&&\Theta(X^g) = \Theta(X)^g \forall g \in G
\eeq
where $G_X$ is the stabiliser of $X$, $G_X = \{g \in G | X^g = X\}$.

Let $T_k$ be the set of sets of size $k$ such that
$P(X)$ is true for all $X \in T_k$, and $T_k$ contains no isomorphs.
The following algorithm generates a set $T_{k+1}$ from $T_k$.

\begin{alltt}   
\(T\sb{k+1}\) := \(\emptyset\)
for \(X in T\sb{k}\) do
  for \(x\) representative in each orbit of \(G\sb{X}\) upon V - X do
    if \(P(X+x)\) and \(x \in \Theta(X+x)\) then
      add \(X+x\) to \(T\sb{k+1}\)
    endif
  enddo
enddo
\end{alltt}

\begin{thm}[Royle\cite{royle98}, McKay \cite{mckay98}]
Let $T_k$ contains exactly one representative from
each $G$--orbit on $k$--sets of $V$ that have property $P$. 
Then the set $T_{k+1}$ generated by the algorithm above 
contains exactly one representative from
each $G$--orbit on $k+1$--sets of $V$ that have property $P$.
\end{thm}

Starting with $T_0 = \{\emptyset\}$ we obtain an orderly algorithm for
constructing one representative of each subset of $V$ with
property $P$.

Royle's technique for orderly algorithmic generation has been applied
to a number of problems, mostly in the area of flocks and related
structures in finite geometries.

One main problem is to define the function $\Theta$.

One property of the  {\tt nauty} package is that it constructs
a \emph{canonical labeling} of a graph. A canonical labeling
uniquely identifies each node of a graph up to automorphisms.

\begin{defn}
A \emph{canonical mapping}  $\alpha$ takes a graph $\Gamma = (N,E)$
and maps 
\begin{eqnarray}
\alpha : {\cal G} &\rightarrow& ({\cal N} \rightarrow \N)\\
\alpha :  \Gamma=(N,E) &\mapsto& \alpha_\Gamma : N \rightarrow\{1,\ldots|N|\}
\end{eqnarray}
s.t.\ for any bijection  $\phi:N\rightarrow M$, $\forall n\in n, m\in M$,
\beq
\alpha_{\Gamma^\phi}(m) = \alpha_\Gamma(n) \Rightarrow
\exists \psi \in Aut(\Gamma), \phi(n^\psi) = m
\eeq
\end{defn}

One such function is available in GAP using the GRAPE package which forms an
interface to the nauty package. With this canonical mapping we can
work wonders; this is in some sense the magic bullet for Royle's algorithm.

\begin{defn}
Let $f:2^V \rightarrow {\cal G}$ be a mapping from subsets of $V$ to
the class of digraphs. Let $X \subseteq V$, $f(X) = (N,E)$
and $f_X:X\rightarrow N$ such that $\forall g \in G$ there
exists an isomorphism $\phi_g:f(X) \rightarrow f(X^g)$ such that
$f_{X^g} \circ g = \phi_g \circ f_X$ and for all $\phi \in Aut(f(X))$
there exists $g \in G_X$ such that $\phi = \phi_g$.
We call such a pair of mappings a \emph{canonical embedding}.
\end{defn}

This definition means that we can embed our sets $X$ into some class
of graphs such that the automorphisms of $X$, $G_X$ and the automorphisms
of the class of graphs coincide properly. We may need to make
a complex embedding in order to prohibit extra automorphisms
from arising.

\begin{prop}
\label{prop_ord_theta}
Let $f$ define a canonical embedding, $\alpha$ a canonical
mapping. Define $\Theta(X) = x^{G_X}$ with 
$\alpha_{f(X)}(f_X(x))$ minimal in $\alpha_{f(X)}(f_X(X))$.
Then $\Theta$ satisfies the requirements for an orderly algorithm.
\end{prop}
\begin{proof}
The first two properties required of $\Theta$ are apparent.

Let $g\in G$ be arbitrary but fixed.
Define $x_1,x_2 \in X$ such that $\alpha_{f(X)}(x_1)$ is minimal and
  $\alpha_{f(X^g)}(x_2^g)$ is minimal. 
Since they
are both minimal in $1\ldots |X|$ they are equal,
so by the canonical mapping property there exists
$\psi \in Aut(f(X))$ such that
\[
\phi_g \psi (f_X(x_1)) =  f_{X^g} (x_2^g)
\]
Since
\[
\phi_g^{-1} f_{X^g} (x_2^g) = \phi_{g^{-1}} f_{X^g} (x_2^g) 
    = f_X {g^{-1}} (x_2^g) = f_X(x_2)
\]
Using the property of a canonical embedding, and the fact that
the automorphism $\psi$ can be represented as $\phi_h$ for 
some $h\in G_X$, we obtain
\[
f_X(x_2) = \phi_h  f_{X}(x_1) = f_X(x_1^h) \Rightarrow x_2 = x_1^h.
\]
We know $\Theta(X) = x_1^{G_X}$ and $G_{X^g} = g^{-1}G_X g$.
Then
\begin{eqnarray}
\Theta(X^g) &=& (x_2^g)^{G_{X^g}} = x_2^{g g^{-1} G_X g} \\
 &=& (x_2^{G_X})^g  = (x_1^{G_X})^g = \Theta(X)^g
\end{eqnarray}
and we are done.
\end{proof}

\begin{defn} A mapping $v:2^V\times V \rightarrow O$, $O$ an ordered set,
 is
a \emph{combinatorial value} if $v(X,x) = v(X^g,x^g)$ 
$\forall g \in G$.
\end{defn}
We usually take $O$ to be the integers, or the set of $n$--tuples
ordered lexicographically.
\begin{corr}
\label{corr_ord_theta}
Let $v_i, i=1,\ldots,n$ be  combinatorial values on $V$.
With terms as in the previous Proposition, requiring
$(v_1(X,x),\ldots,v_n(X,x),\alpha_{f(X)}(f_X(x))$ to be minimal
in $\{(v_1(X,x),\ldots,v_n(X,x),\alpha_{f(X)}(f_X(x)) : x \in X\}$
gives a $\Theta$ that satisfies the requirements of an orderly algorithm.
\end{corr}
\begin{proof}
If $x$ is such that $v_i(X,x)$ is uniquely minimal, then we have
determined a node uniquely. Thus $\Theta(X) = \{x\}$,
so $\Theta(X^g) = \{x^g\} = \Theta(X)^g$ and we are done.

Otherwise we select our minimal $x$ as in Proposition \ref{prop_ord_theta}
from a union of orbits of $G_X$ on $X$. This union of orbits is covariant
across isomorphisms, so we can use the same argument as in
the proof of Proposition \ref{prop_ord_theta}.
\end{proof}

We can compute this $n+1$--tuple stepwise; if $v_i(X,x)$ is
not minimal then we know that $x$ is not in the orbit $\Theta(X)$.
If  $v_i(X,x)$  is uniquely minimal for all $x\in X$, 
then it is in  $\Theta(X)$ and we can stop
and accept the extension (because $x \in \Theta(X)$), without having to compute the 
following $v_j(X,x)$ or the canonical labeling.
Otherwise we carry on to the next element of the tuple.
Only in some cases do we need to compute the canonical labeling.

We could compute the orbits of $G_X$ on the elements $x$
with $(v_1(X,x),\ldots,v_n(X,x))$ minimal. If there is only one
orbit, then $\Theta(X)$ must be this orbit. This removes further
determinations of the canonical mapping. However, in GRAPE we
obtain the automorphism group of the graph and the canonical
mapping at the same time, so we do not save any computation.

\section{Central Groupoids}

Evans  \cite{evans67a}
 investigates the various products that can
be defined on the set $S = A \times A$ for some
set $A$.
A well--known example of such a construction is
the rectangular semigroup \cite{mclean54}, where one takes
a pair of sets $A,B$, forms the product $S=A \times B$
and the product $(a,b) * (c,d) = (a,d)$.
This is an idempotent 
semigroup satisfying the axiom $a*b*a=a$ for all $a,b \in S$,
and every semigroup satisfying this equation is of this form.

Evans looked at all the possible products on $S=A\times A$,
and found that other than the rectangular semigroups, the
only other interesting examples were defined by
\begin{eqnarray}
(a,b) \bullet (c,d) = (b,c).
\end{eqnarray}
This operation satisfies the identity
$(a\bullet b)\bullet (b\bullet c)=b$.

\begin{defn} 
A {\em Central Groupoid}
 is a  $(2)$-algebra $(S,\bullet)$ 
satisfying the   axiom:
\begin{eqnarray}
  (a\bullet b)\bullet (b \bullet c) = b
\end{eqnarray}
\end{defn}

The examples used by Evans are referred to as the \emph{natural} central groupoids.
All natural central groupoids have, by construction,
order equal to a square.
In a natural central groupoid, we  see
that $(a,b) \bullet (a,b) = (a,b)$ iff $a=b$,
thus we have $|A|$ idempotents
in a natural central groupoid of order $|A|^2$.
These are general results for central groupoids.

\begin{thm}[Evans, Knuth]
If $(S,\circ)$ is a finite central groupoid, then
$\vert S \vert = n^2$ for some integer n. 
For every positive integer $n$ there exists a central groupoid
of order $n^2$. In any finite central groupoid of
order $n^2$, the number of idempotents is $n$.
\end{thm}

Note that
the first result appears as  Corollary \ref{rs_cor_evans}, 
 the second follows from the example above.
The third result is very difficult to show in an algebraic
setting, we need to move over to a matrix theoretic 
setting in order to prove it \cite{knuth70}.

In \cite{knuth70}, Knuth investigated
various aspects of central groupoids.
Most importantly, Knuth found two  models
for central groupoids, one being a digraph model,
the other being a model based upon the 
$\{0,1\}$--matrices that are the incidence matrices of these
digraphs.
The incidence matrices made an interesting contribution
to a question posed by Hoffman in \cite{hoffman67}: 
which $\{0,1\}$--matrices $A$ have the property that
$A^2=J$, where $J$ is the matrix consisting entirely of
ones?

In \cite{shader74} Shader demonstrates that
 non--natural
central groupoids exists for all orders $n^2$, 
$n \geq 3$.
The question of an exhaustive list of central groupoids,
or equivalently an exhaustive list of matrices
$A$ with $A^2=J$, remains open.

A digraph has the \emph{unique path property} of length $n$
$UPP_n$, if there is a unique path of length $n$ between
any two nodes \cite{fiol85,mendelsohn70}. 
For $n=1$ we obtain complete graphs.
If $(S,\cdot)$ is a central groupoid, then the digraph
with node set $S$ and edges $\{(a,a\cdot b): a,b \in S\}$ is a
$UPP_2$ graph.
If $(N,E)$ is a $UPP_2$ graph, then for every pair of nodes
$a,b \in V$, there is a unique node $c$ that is the midpoint
of the unique path from $a$ to $b$. Defining $a\circ b := c$
gives us a binary operation and $(V,\circ)$ is a central groupoid.
Thus the determination of
$UPP_2$ graphs and central groupoids is equivalent. If $A$ is the
incidence matrix of a $UPP_2$ graph, then $A^2=J$ 
where $J$ is the matrix consisting entirely of
ones, so we have the solutions to Hoffman's problem.

Many investigations of this and related problems have been made,
in particular looking at special classes of examples from 
circulant matrices
\cite{lam75,lam77,lamlint78,ma84,ryser70,wang80,wang82,wulihuang99,wu02}.

\section{Direct Implementation}

This problem of determining all $UPP_2$ graphs seemed like
an easy problem for an orderly algorithm approach. We proceeded as follows.

Let $V$ be the set of directed edges in a graph of order $n^2$.
Let $G$ be the symmetric group on $n^2$ points acting upon these
edges naturally. Let $P(X)$ be the property of a graph defined by the edges
$X$ having at most one path of length 2 between any pair of nodes.
$\Theta(X)$ can be directly defined using the canonical labeling property
described above on the graph induced by $X$.

We can remove many possible extensions using other properties of 
central groupoids; for instance we know (see below) that the
in and out valency of each node in a $UPP_2$ graph is $n$. Thus
we do not add edges between nodes that have already reached that limit.

Although the algorithm, when run on the problem for $n=2$, terminates
quickly with the correct answer, attempts (with lots of clever implementation
tricks as suggested by Royle) for $n=3$ ran for several weeks
before we terminated them. 

The problem seems to be that the number of stages in the algorithm
is $n^3$. In the process of building partial solutions, many
partial examples are constructed that cannot be extended to a full
solution. It seems that the recursion depth $n^3$ allows a
lot of extraneous and unproductive branching, bogging the
process down and preventing the algorithm from working effectively.

Thus we resorted to a different approach. The structures that we
introduce in the next section are generalisations of central groupoids,
yet have a combinatorial structure that is somewhat simpler.

\section{A Generalisation}

We introduce a generalisation of the central groupoid structure.
Not all results will be proved here, proofs can be found in 
\cite{boykett04}.

\begin{defn} 
\label{axnoq}\index{Semicentral Bigroupoid}
A {\em Semicentral Bigroupoid} is a  $(2,2)$-algebra $(S,\bullet,\circ)$ 
satisfying the following axioms:
\begin{eqnarray}
(a\bullet b)\circ (b \bullet c) = b\\
(a\circ b)\bullet (b \circ c) = b 
\end{eqnarray}
\end{defn}

If $(S,\bullet)$ is a central groupoid, then $(S,\bullet,\bullet)$
is a semicentral bigroupoid, so this is a proper generalisation.

Note also that the definition is completely symmetric
in $\bullet$ and $\circ$, i.e.\  $(S,\circ,\bullet)$ is a 
semicentral bigroupoid iff $(S,\bullet,\circ)$ is.
The \emph{dual} of a   \scbg\ $(S,\bullet,\circ)$
is $(S,\circ,\bullet)$.
It is often not necessary to prove results
for both operations, as they carry across by duality.

\begin{ex}
\label{ex_assoc}
Let $A$,$B$ be two sets, and let $Q=A\times B$. Define
\begin{eqnarray}
(a_1,b_1) \bullet (a_2,b_2) &=& (a_1,b_2) \\
(a_1,b_1) \circ (a_2,b_2) &=& (a_2,b_1) 
\end{eqnarray}
Then $(Q,\bullet,\circ)$ is a \scbg. 
\end{ex}

In the above,  $(Q,\bullet)$
defines   a rectangular semigroup. 
This  corresponds
to the natural central groupoids, in that it can be constructed
as a ``product of points'' \cite{evans67a}, see also Lemma 
\ref{lemma_nat_lift}.
Below we will see that
it is the only associative \scbg.

\begin{lemma}
\label{lemma_anticomm}
\label{lemma_idem}
Let $(S,\bullet,\circ)$ be a \scbg.
Both the operations are anti-commutative, that
is, $a \bullet b = b \bullet a \Rightarrow a = b$ and
similarly for $\circ$.
Also $ a \bullet a = a$ iff $a\circ a = a$,
thus $(S,\bullet)$ is idempotent iff $(S,\circ)$ is idempotent.
\end{lemma}

These results  follow by direct calculation. Note that anticommutative
is a stronger condition than noncommutative. In the following,
we will often omit $\bullet$ and represent the operation by juxtaposition
where no confusion would result.

We can take any semicentral bigroupoid and ``bend''
it to get another semicentral
bigroupoid.
This method can be used to find new examples of
semicentral
bigroupoids.

\begin{prop}
\label{propperm}
If  $(S,\bullet,\circ)$ is a semicentral bigroupoid, 
and $\phi : S \rightarrow S$ is a permutation of $S$, then
the algebra $(S, *, +)$ with
\begin{eqnarray} 
a * b = \phi^{-1}(a\bullet b)\\
a + b = \phi (a) \circ \phi(b)
\end{eqnarray}
is also a semicentral bigroupoid.
\end{prop}

The calculation behind this result is mechanical. 
Note that the new semicentral bigroupoid will 
in general not be isomorphic to the old one; see  
 Proposition \ref{prop_iso}.

\begin{defn}
The {\em lifting of $(S,\bullet,\circ)$ by $\phi$} is the
algebra $(S,*,+)$ defined above.
The {\em square map} $\phi_\bullet$ of $(S,\bullet,\circ)$
is $\phi_\bullet: x \mapsto x\bullet x$
\end{defn}

The square map, $\phi_\bullet : a \mapsto a\bullet a$ is a permutation:
\begin{eqnarray}
&&\phi_\bullet(a) = \phi_\bullet(b)
  \Rightarrow a \bullet a = b\bullet b \nonumber \\
&&\;\;\;\;\Rightarrow (a \bullet a) \circ (a \bullet a) = 
		(b\bullet b) \circ (b \bullet b) 
  \Rightarrow a = b
\end{eqnarray}
If we lift by the square map  then
the derived operation $*$ is idempotent:
\begin{eqnarray}
a*a = \phi_\bullet^{-1} (a \bullet a ) = \phi_\bullet^{-1}\phi_\bullet(a) = a.
\end{eqnarray}
This will be referred to as the {\em idempotent lifting}
of a semicentral bigroupoid. 
If we then lift the resulting idempotent \scbg\ by
the permutation $\phi_\bullet^{-1}$ then we will obtain $(S,\bullet,\circ)$.

It is clear that the lifting operation
is invertible  and that
every  \scbg\ is the lifting of an idempotent \scbg\
by a permutation that becomes the inverse of
the square map in the lifting.

That we can have any permutation as the square map
in a \scbg, and thus any number of idempotents,
 is a contrast to the case for a 
central groupoid, where there are
exactly $\surd\vert S \vert$ idempotents
for $S$ finite.

We see the following.

\begin{prop}
Every \scbg\ $(S,\bullet,\circ)$ can be uniquely represented as 
an idempotent \scbg\ and a permutation in $Symm(S)$. Conversely, every
such pair gives a \scbg\ that is idempotent iff
the permutation is trivial.
\end{prop}

Every central groupoid is the lifting of an idempotent \scbg\
such that the two operations in the lifting are identical.
In Proposition \ref{prop_iso} we will see exactly when 
two \scbgs\ are isomorphic, 
based upon the isomorphism of their idempotent representatives and 
relations between their square maps.
We will focus upon idempotent \scbgs\ for now.

Let's apply this to a central groupoid.
As mentioned above, if $(S,\bullet)$ is
a central groupoid, then  $(S,\bullet,\bullet)$ is
a semicentral bigroupoid. 
Take the example of the natural central 
groupoid of order 4. The elements are $\{aa,ab,ba,bb\}$ with operation
$x_1x_2 \bullet x_3x_4 = x_2x_3$.
The square map is the permutation
$\phi_\bullet = (ab\; ba)$.
This is the permutation that reverses the entries in the
product, i.e.\  $\phi_\bullet: xy \mapsto yx$.
If we construct the multiplication tables for the lifting by $\phi_\bullet$,
then we obtain the following:
\begin{eqnarray}
\begin{tabular}{c | c c c c}
* & aa & ab & ba & bb\\
\hline
aa & aa & aa & ba & ba \\
ab & ab & ab & bb & bb \\
ba & aa & aa & ba & ba \\
bb & ab & ab & bb & bb \\
\end{tabular}
\;\;\;\;\;
\begin{tabular}{c | c c c c}
+ & aa & ab & ba & bb\\
\hline
aa & aa & ab & aa & ab \\
ab & aa & ab & aa & ab \\
ba & ba & bb & ba & bb \\
bb & ba & bb & ba & bb 
\end{tabular}
\label{eq_ncg4_lift}
\end{eqnarray}
We see that this lifting is  an idempotent semicentral bigroupoid,
in fact an associative one.

In general one can make the following statement.
\begin{lemma}
\label{lemma_nat_lift}
The idempotent lifting of a natural central groupoid
is an associative semicentral bigroupoid.
\end{lemma}
The proof is mechanical.

In Lemma \ref{lemma_anticomm} above, we saw that the
operations in a  \scbg\ are anticommutative.
In \cite{mclean54}, McLean shows that the only
anticommutative semigroups are the rectangular 
semigroups. Thus if $(S,*,+)$ is a \scbg\ with $(S,*)$ associative,
then $(S,*)$ is a rectangular semigroup, as is $(S,+)$, with a structure
as in Example \ref{ex_assoc}.

\begin{defn}
A {\em Rectangular Structure} 
on a set $S$, called the {\em base set}, 
is a collection $\cal R$ of ordered pairs of subsets, 
called {\em rectangles},
of $S$, such that
\begin{eqnarray}
&\forall (s,t)\in S^2 \;\exists !\; R \in {\cal R} \mbox{ such that }
	(s,t) \in R \label{rs_axiom1}\\
&\forall R,Q \in {\cal R}, \vert R_1 \cap Q_2 \vert = 1 \label{rs_axiom2}
\end{eqnarray}
where we identify $R = (R_1,R_2) = R_1 \times R_2$.

We say two rectangular structures are {\em isomorphic} if there is
an invertible map between the base sets that preserves
rectangles.
\end{defn}

As an example, take two sets $A,B$. Define $S = A\times B$, and
for all $(a,b) \in S$ define 
$R_{(a,b)} = ( \{a\}\times B, A \times \{b\})$.
Then 
\begin{eqnarray}
 {\cal R} = \{ R_{(a,b)} \vert a \in A, b \in B \}
\end{eqnarray}
is a rectangular structure on $S$.
For any $((a,b),(c,d)) \in S^2$, $((a,b),(c,d)) \in R_{(a,d)}$,
and this rectangle is  unique.
Let $R = R_{(a,b)}$, $Q = R_{(c,d)}$ be two rectangles in $\cal R$.
Then $R_1 = \{a\}\times B$ and $Q_2 = A \times \{d\}$,
so $|R_1 \cap Q_2| = |\{\{a\}\times\{d\}\}| = 1$.

Define
\begin{eqnarray}
\rho : S &\rightarrow & {\cal P} (S^2) \\
	x &\mapsto & \{(ax,xb) | a,b \in S\}
\end{eqnarray}
Where ${\cal P}(X)$ denotes the power set of $X$.

\begin{lemma}
\label{lemmrhopart}
\label{lemma_scbg_john}
$\rho(S) = \{\rho(s) \vert s \in S \}$ is a partition of $S^2$.
For every $x \in S$ there exists $A,B \subseteq S$ such that
\begin{itemize}
\item $\rho(x) = A \times B$.
\item $\vert A \cap B\vert = 1$.
\item $B \circ A = S$.
\item $A \circ B =  \{x\}$
\end{itemize}
Thus ${\cal R} = \{\rho(s):s\in S\}$ is a rectangular structure.
\end{lemma}
\begin{proof}
The proof of the 4 points  is mechanical calculation 
with $A=S\bullet x,\ B=x\bullet S$.

For any $a,b \in S$, $(a,b) \in \rho(a\circ b)$ so  (\ref{rs_axiom1})
is satisfied.

Then for any pair of rectangles  $Q,R \in {\cal R}$, 
note that there are some
$x,y \in S$ such that $Q_1 = S\bullet x, R_2 = y \bullet S$.
If $a \in Q_1 \cap R_2$, then $a = b\bullet x = y \bullet c$
for some $b,c \in S$. Then by
Corollary \ref{corrswap} above, $a=y \bullet x$, that is,
$Q_1 \cap R_2 = \{y\bullet x\}$, so condition  (\ref{rs_axiom2})
is satisfied.
\end{proof}

\begin{corr}
\label{corrswap}
Let $(S,\bullet,\circ)$ be a \scbg.
If $a\circ b = c \circ d = x$, then $a \circ d = c \circ b = x$,
and similarly for $\bullet$.
\end{corr}
\begin{proof}
By Lemma \ref{lemma_scbg_john} above, $\rho(x) = A \times B$,
$A = S\bullet x$, $B = x\bullet S$.
Then
\begin{equation}
a = (a\circ a)\bullet(a \circ b) = (a\circ a)\bullet x \in Sx
\end{equation}
similarly $c \in S\bullet x$, $b,d \in x\bullet S$.
Then $A\circ B = (S\bullet x) \circ (x\bullet S) =\{x\}$, so
\begin{eqnarray}
a \circ d \in A \circ B &\Rightarrow & a\circ d = x \\
c \circ b \in A \circ B &\Rightarrow & c \circ b = x
\end{eqnarray}
\end{proof}

Such a property is seen also in the L--groupoids
in \cite{guggenbergerrung03}.

The \emph{format} of a rectangle $(A,B)$ is the 
ordered pair $(\vert A \vert,\vert B \vert)$.

Define the map:
\begin{eqnarray}
 d: \cal R &\rightarrow& S \label{d_defn}\\
    R &\mapsto& r \mbox{ where } \{r\} = R_1 \cap R_2
\end{eqnarray}
This map is well defined since for every $R \in {\cal R}$,
$|R_1 \cap R_2| = 1$ by (\ref{rs_axiom2})  above,
so $R_1 \cap R_2 = \{r\}$ for some unique $r$.

By (\ref{rs_axiom1}) above, $(r,r)$ is in a unique
rectangle, so this map is bijective 
and $\vert {\cal R} \vert = \vert S \vert$.

\begin{prop}[\cite{boykett04} Proposition 8]
\label{prop_rs_size}
If $\cal R$ is a rectangular structure with base set $S$, and
$R = (R_1,R_2) \in \cal R$ is some rectangle, then
$\vert R_1\vert \vert R_2 \vert = \vert S \vert = \vert \cal R \vert$.
Moreover, for any other rectangle $Q = (Q_1,Q_2) \in \cal R$,
$\vert R_1 \vert = \vert Q_1 \vert$, 
i.e.\ all rectangles have the same format.
\end{prop}

A \rs\ ${\cal R} = \{A_i\times B_i \}$ is \emph{right (left) partitioned} 
if the sets $B_i$ ($A_i$) form a partition of $S$. 
If it is partitioned on both sides, we call it doubly partitioned.
If $(S,*,+)$ is a \scbg\ then $\cal R ^*$ is left partitioned
iff $\cal R^+$ is right partitioned (and vice versa).
Using partitions and orthogonal partitions, we can construct
many well-behaved \rss. This will however not concern us, as
we will see later in the comments before Lemma \ref{lemma_nat_assoc}.

From  a rectangular structure $\cal R$, 
using the bijection $d$ from equation (\ref{d_defn}) above
and denoting by $R(s,t)$ the unique rectangle on the pair $(s,t)$
guaranteed by (\ref{rs_axiom1}),
 define
\begin{eqnarray} 
 \bullet : S \times S &\rightarrow& S \label{eq_algg1}\\
  (s,t) &\mapsto& u \mbox{ where } \{u\} = (d^{-1}(s))_2 \cap  (d^{-1}(t))_1  \nonumber\\
 \circ : S \times S &\rightarrow& S \label{eq_algg2}\\
  (s,t) &\mapsto& d(R(s,t)) \nonumber
\end{eqnarray}
as binary operations on $S$.

\begin{prop}[\cite{boykett04} Proposition 9, Lemma 5]
\label{prop_rect_to_scbg}
The algebra $(S,\bullet,\circ)$, with operations defined
as in (\ref{eq_algg1}),(\ref{eq_algg2}) above, is an idempotent semicentral bigroupoid. The idempotent \scbg\ is associative iff the \rs\ is partitioned on
both sides.
\end{prop}

If we call the format of an operation table the format of
the derived rectangular structure, we get the following.

\begin{corr}
If $(S,\bullet,\circ)$ is a semicentral bigroupoid with
 format $(a,b)$ for the $\bullet$ operation table,
then the format of the $\circ$ operation table is $(b,a)$.
\end{corr}
\begin{proof}
Let $(c,d)$ be the format of the $\circ$ operation table.
For any $x \in S$, $\rho(x) = S \bullet x \times x \bullet S$
is the rectangle filled with $x$ in the $\circ$ operation table.
Thus $d = \vert x \bullet S \vert$, so there are 
$d$ rectangles in the $x$ row of the $\bullet$
operation table, all of which have the same format $(a,b)$.
Thus $db = \vert S \vert$. But $\vert S \vert = ab$ so $a=d$.
Similarly $b=c$, i.e.\ the format of the $\circ$ operation
table is $(b,a)$.
\end{proof}

\begin{corr}[\cite{evans67a} Theorem 1]
\label{rs_cor_evans}
A finite central groupoid $(S,\bullet)$ has square
order.
\end{corr}
\begin{proof}
If $(S,\bullet)$ is a central groupoid, then $(S,\bullet,\bullet)$
is a \scbg. Thus the formats of the
operations are $(a,b)$ and $(b,a)$, but these are identical,
so $a=b$ and $\vert S \vert = ab = a^2$.
\end{proof}

We now investigate the isomorphism of \scbgs.

\begin{prop}[\cite{boykett04} Lemma 16, Proposition 11]
\label{prop_bertl}
Every  semicentral bigroupoid $(S,\bullet,\circ)$   has an 
associated rectangular structure that is constant across liftings.
Two idempotent  semicentral bigroupoids are isomorphic iff
the associated rectangular structures are isomorphic.
Isomorphic (nonidempotent) \scbgs\ have isomorphic \rss.
\end{prop}

The construction of a rectangular structure from a  
semicentral bigroupoid above and the
reverse construction in Proposition \ref{prop_rect_to_scbg}  are
complementary in that given an idempotent  semicentral bigroupoid
$(S,*,+)$, the  semicentral bigroupoid $(S,\bullet,\circ)$ derived from the
associated rectangular structure is the same, 
i.e.\ $a*b = a \bullet b,\; a + b = a \circ b$.

 Rectangular structures and idempotent  semicentral bigroupoids are
essentially equivalent. 
Now to look at a similar result for non--idempotent  semicentral bigroupoids.

Let $Symm_{RS}(S)$ be
the symmetry group of the rectangular structure of $(S,\bullet,\circ)$.

\begin{prop}[\cite{boykett04} Propositions 12 and 13]
\label{propcon}
\label{prop_iso}
Two  semicentral bigroupoids 
$(S,\bullet,\circ)$ and $(T,*,+)$ are isomorphic, 
with isomorphism $\beta:S \rightarrow T$, iff
their idempotent liftings are isomorphic by $\beta$, and
$\beta \phi_\bullet = \phi_*\beta$ for the
square maps $\phi_\bullet$ and $\phi_*$.
Two liftings of an idempotent  semicentral bigroupoid by  $\phi, \bar\phi$ are
isomorphic
iff
the $\phi, \bar\phi$ are conjugate by an element from $Symm_{RS}(S)$.
\end{prop}

Thus we can determine exactly when two \scbgs\ are isomorphic:
they must have isomorphic \rss\ and the liftings must be
conjugate by the symmetric group of the \rs.
We will see later how to use this for the determination
of central groupoids.

We generalise the $UPP_2$ condition as follows.
Take a fixed set of vertices, and look at two directed graphs
on this set, $G_R$ and $G_B$.
Call these as the red and blue graphs
respectively.
The problem is to arrange these graphs such that,
when we superimpose them, obtaining a multigraph, 
there is a unique directed path of length
2 coloured blue--red between any two nodes, and a unique
directed path coloured red--blue.
We refer to these as
\emph{symmetric unique 2--coloured 2--path graphs}.
We must allow (in general) that two edges of differing colours
exist between the same vertex pair. Thus we should speak
of multigraphs.
Thus the class of multigraphs might be termed
\emph{symmetric 2--coloured unique 2--path multigraphs}.

We note that the incidence matrices of these graphs have the property that
$AB = BA = J$ where $A,B$ are the incidence matrices of the two graphs
and $J$ is the matrix consisting of all 1s.
These are generalisations of the matrices introduced
by Hoffman in \cite{hoffman67}.

\begin{ex}[Construction]
\label{ex_constr}
Take a semicentral bigroupoid, define coloured edges
$a \rightarrow_{blue} b$ and 
$a \rightarrow_{red} b$ by
\beq
a \rightarrow_{red} b & \Leftrightarrow &
	a \bullet c = b \;\exists c \\
a \rightarrow_{blue} b & \Leftrightarrow &
	a \circ c = b \;\exists c 
\eeq

Conversely, given a graph pair define
a \scbg\ as follows. 
Let $S$ be the vertex set of the graphs.
For any $a,b \in S$, let $c_r,\; c_b$ be the unique
vertices on the paths between $a$ and $b$:
\beq
a \rightarrow_{red} c_r \rightarrow_{blue} b \\
a \rightarrow_{blue} c_b \rightarrow_{red} b 
\eeq
Define
\beq
a\bullet b &=& c_r \\
a\circ b &=& c_b
\eeq
and we have an $(2,2)$--algebra $(S,\bullet, \circ)$.
\end{ex}

\begin{prop}
If $(S,\bullet,\circ)$ is a \scbg, then the construction above
determines a symmetric 2--coloured unique 2--path multigraph.
Similarly, given a symmetric 2--coloured unique 2--path multigraph,
the construction above defines a \scbg. The constructions are inverses
of one another.
\end{prop}
The proof is mechanical.

Note that if some element $a$ is idempotent,
then $a\bullet a = a$ so there is a red loop
arc on the node $a$, similarly a blue loop arc.
Thus if $S$ is idempotent, every node is
a loop node.

Consider two categories, ${\cal S}$ of semicentral bigroupoids,
and $\cal G$ of symmetric 2--coloured unique 2--path multigraphs.
Consider the functor from $\cal S$ to $\cal G$
as described in Example \ref{ex_constr} 
above, and take an exact sequence
$A \rightarrow_f B \rightarrow_g C$ in $\cal S$.
Since the mappings $f,g$ operate on elements of the
semicentral bigroupoids, they  operate
on vertices of the graphs under the functor,
carrying arcs across with them.
Thus $range(f)= ker(g)$ in $\cal G$, so
the functor is exact.
Thus we get the result:

\begin{lemma}
\label{lemma_exact_functor}
Isomorphic semicentral bigroupoids give isomorphic
multigraph pairs by the   construction in Example \ref{ex_constr} 
above.
\end{lemma}

Thus the correspondence between
multigraph pairs and \scbgs\ is an equivalence. 
The following  shows that
the automorphism group of a \scbg\ can be obtained using
the automorphism groups of the associated graphs.

\begin{lemma}[\cite{boykett04} Lemma 14]
\label{lemma_autgp_graphs}
If $G_b,G_r$ are the graphs defined by a  semicentral bigroupoid $S$, then
\begin{equation}
Aut(S) = Aut(G_b) \cap Aut(G_r)
\end{equation}
\end{lemma}

If we look at idempotent \scbgs, then the associated multigraph pairs
have loop edges of both colours on each node. However they
do not have any other repeated edges. We can then ignore the
loop edges for determining the automorphism groups and other
such properties.

\section{Applications to Central Groupoids}

We see that the \scbgs\ have a lot of structure. This section investigates
what this means for
central groupoids. 
We have already seen that the format of a 
central groupoid must be $n \times n$. A central groupoid is a 
lifting of a square idempotent \scbg\ by a
permutation. 

\begin{lemma}
Let $(S,+,*)$ be an idempotent \scbg, $\phi$ a permutation of $S$.
Then  
the lifting of $S$ by $\phi$ is
a central groupoid iff
$\phi:(S,+)\rightarrow (S,*)$
is an isomorphism of order 2.
\end{lemma}
\begin{proof}
($\Leftarrow$)
Let $(S,\bullet,\circ)$ be the lifting of  $(S,+,*)$ by $\phi$.
We know that
\begin{eqnarray*}
a \bullet b &=& \phi(a+b) \\
a \circ b &=& \phi a * \phi b
\end{eqnarray*}
By the isomorphism property of $\phi$ then
$a \bullet b = \phi(a+b) = \phi a * \phi b = a \circ b$
so $(S,\bullet,\bullet)$ is a \scbg\ so $(S,\bullet)$ is a central groupoid.

($\Rightarrow$)
Let $(S,\bullet,\circ)$ be the lifting of  $(S,+,*)$ by $\phi$,
$(S,\bullet)$ = $(S,\circ)$ a central groupoid.
Then $a\bullet a = \phi(a+a) = \phi(a)$, 
$a=(a\bullet a)\bullet (a\bullet a)
= \phi \phi (a)$ for all $a\in S$ so $\phi$ is of order 2.
By the central groupoid property, $\bullet = \circ$, so
\[ \phi(a+b) = a\bullet b =a \circ b = \phi(a) * \phi(b)\]
so $\phi$ is an isomorphism.
\end{proof}

The two graphs defined above from a \scbg\ are defined by the
operations $+$ and $*$, so this means that the two graphs are
isomorphic with isomorphism $\phi$.

If the \rs\ underlying a central groupoid is partitioned, then
the partitions are apparent on one side of the multiplication table
of $*$ but on the other side of the multiplication table of $+$.
However these operations are isomorphic, so if the \rs\ is
partitioned, then it must be partitioned on both sides.

\begin{lemma}
\label{lemma_nat_assoc}
If $S$ is a square associative \scbg\ then there is only
one central groupoid lifting of it and it is the
natural central groupoid.
\end{lemma}
\begin{proof}
The automorphisms of a rectangular semigroup $S=A \times B$
are direct products of permutations of $A$ and $B$. 

Let $(S,+,*)$ be an associative \scbg\ of format $n\times n$,
$S = A\times A$.
Let $\phi$ be an isomorphism of $(S,+)$ onto $(S,*)$. 
The map $\sigma:(a,b)\mapsto (b,a)$ is also such an isomorphism,
thus the composition $\sigma\phi$ is an automorphism
of the rectangular semigroup $(S,*)$. Thus $\phi(a,b) = (\phi_1b,\phi_2a)$
for some permutations $\phi_1,\phi_2$ of $A$.

Let $\phi$ and $\psi$ be two permutations of $S$ giving
central groupoid liftings.
We know that $\phi$ and $\psi$ have $n$ fixed points and $\frac{n^2-2}{2}$
swaps, so all of $\phi_1,\phi_2,\psi_1,\psi_2$ must have the same structure.
Thus
there are permutations $\alpha_1,\alpha_2$ of $A$ such that
$\phi_1^{\alpha_1} = \psi_1$ and $\phi_2^{\alpha_2} = \psi_2$.
Since $(\alpha_1,\alpha_2)$ is an isomorphism of $S$,
we know that the liftings by $\phi$ and $\psi$ are
isomorphic using Proposition \ref{propcon}.

It is clear that the lifting of $S$ by $\sigma$ is 
the natural central groupoid.
\end{proof}

Given a square \rs\ $\mathcal R$, we can determine whether a 
central groupoid can
be defined from it as follows:
\begin{itemize}
\item if $\mathcal R$ is partitioned, it must be partitioned on both sides.
Which means it must be the \rs\ made by products of points, and the
only central groupoid that can be obtained from it is the natural one.
\item The graph pair generated from $\mathcal R$ must be isomorphic.
\item The set of isomorphisms can then be reduced to those of
order 2. These are the candidate liftings for the central groupoid.
\item These candidates are then separated into orbits under the conjugation
operation by the automorphism group of the \rs. Each orbit 
representative gives
us one central groupoid, and all such central groupoids are non-isomorphic.
\end{itemize}

\section{Generating and Filtering Process}

In this section we outline the techniques that we used to
obtain exhaustive enumerations of \rss\ and to then
filter out the appropriate examples for the construction
of central groupoids.

\begin{defn} For positive integers $n,m$, an {\em $n\times m$ rectangle} 
is a pair of sets $(R_1,R_2)$
with $\vert R_1 \vert = n, \vert R_2 \vert = m$,
$R_i \subset \{1,\ldots,(nm)\},\; \vert R_1 \cap R_2 \vert = 1$
\end{defn}

For instance $(\{1,4\},\{1,2,3,5\}),\;(\{1,2\},\{2,3,4,5\})$ 
are examples of $2\times 4$ rectangles.

\begin{defn}
A {\em $n \times m$ Partial Rectangular Structure} 
$P$ is a collection of $n\times m$
rectangles such that
\begin{itemize}
\item For all $Q,R \in P, \vert Q_1 \cap R_2\vert = 1$
\item For all $a,b \in \{1,\ldots,nm\}$ there is at most one $R \in P$ such
that $ a \in R_1, b \in R_2$.
\end{itemize}
\end{defn}

One can easily see that a rectangular structure is a partial rectangular structure.
Partial rectangular structures  generalise  rectangular
structures with regard to the 
covering requirement (equation (\ref{rs_axiom1})).
A \emph{full} partial rectangular structure is one with 
$nm$ rectangles. This is a rectangular structure.

We can obtain a graph pair from a partial rectangular structure
in the same way as we obtain it from a rectangular structure.

 We can combine the two graphs
into one graph, the isomorphism of two graph pairs being equivalent to the
isomorphism of two combined graphs. 

A direct method is to generate a graph from the graph pair by
having a coloured graph on $n$ nodes. We colour the
edges depending upon which graph they come from.
We do not get repeated edges, so this is well defined.
However GRAPE does not allow us to pass colourings to the
graph automorphism test, so we have to use a more intricate method.

If the graph pair is $A$, $B$ on node set $\{1,\ldots,n\}$, 
then using the labels
$a$ and $b$ we construct a graph with nodes
\beq
\{(a,x),(b,x)| x \in \{1,\ldots,n\}\}
\eeq 
and edges
\beq
\{((a,x),(a,y))| (x,y) \in Edges(A)\} \cup \nonumber \\
\{((b,x),(b,y))| (x,y) \in Edges(B)\} \cup \nonumber \\
\{((a,x),(b,x))| x \in\{1,\ldots,n\}\} \} \cup \nonumber \\
\{((a,x),(a,x))| x \in\{1,\ldots,n\}\} \}
\eeq 
It is easy to see that the graphs  constructed from 
two graph pairs are isomorphic iff the graph pairs
are isomorphic.
More importantly, the canonical labeling property of nauty, accessed via GRAPE, can
be used on this graph. This is a canonical embedding, as defined above.

Our situation is as follows. The set $V$ is the set of
$n \times m$ rectangles. $P$ is the property of being 
a partial rectangular structure. We want to find the
set $T_{nm}$ of full partial rectangular structures,
i.e.\ \rss.

Note that we assemble only rectangles extending a 
given PRS that will satisfy $P$. We never construct the complete set of
all rectangles, only those that are possible extensions, with
a given middle.

We define $\Theta$ as in Corollary \ref{corr_ord_theta}.
We label rectangles $R$ by their \emph{middle}, the unique element 
$m(R) \in R_1 \cap R_2$.  Then the map $F_{\cal R}(R)$ 
where ${\cal R}$ is a partial \rs\ and
$R\in {\cal R}$, takes $R$ to the node labeled with $(a,m(R))$ 
in the graph defined above.

Initially we used the combinatorial values $v_1({\cal R},R)
= |\{Q\in {\cal R}: m(R) \in Q_1\}|$ and $v_2({\cal R},R)
= |\{Q\in {\cal R}: m(R) \in Q_2\}|$.
The addition of a measure of 
intersection numbers, the ordered pair of multisets 
\[
v_3({\cal R},R)
=(\{|R_1 \cap Q_1| : Q \in {\cal R}\},\{|R_2 \cap Q_2| : Q \in {\cal R}\})
\]
reduced computation time from 6 hours to 39 minutes.

\section{Results}

We have run the algorithms for examples up to order 12. Here 
we are only interested in the case of orders 4 and 9. For order 4 we obtain 
three \rss, one doubly partitioned and two partitioned. Thus only
the one (natural) central groupoid can be obtained.
The order 9 search ran in a little over 39 minutes, drastically faster than
the algorithm we attempted earlier. There are 184 \rss\ of order 9.

We proceeded as follows. We filtered the partitioned
\rss\ out, leaving 67 examples. We then constructed the
graph pairs and found only 7 pairs with isomorphic graphs.
Of these, three had no isomorphisms of order 2, so they could be immediately
excluded.
For a further two examples, we found a unique isomorphism of order
2. 

The remaining two \rss\ both had two isomorphisms of order 2.
Thus we obtained the automorphism group of the \rs\ and determined
the orbits of these permutations under the group of automorphisms
by conjugation, as described in Proposition \ref{propcon} above. 
For one example,
both permutations were in the same orbit, so they give isomorphic
liftings. In the other example, both permutations were fixed by the
automorphism group. Thus they lead to non-isomorphic liftings, 
i.e.\ non-isomorphic
central groupoids.

The doubly partitioned \rs\ has only one central groupoid 
lifting, the natural one.
Thus we obtained a total of one natural and 5 unnatural 
central groupoids of order 9, in accordance with the 
results in \cite{curtis04}.

The numbers in the following refer to the indices in the 
list of rectangular structures obtained. The data is available from the author.
The last two examples can be seen to arise from the
same rectangular structure. Note that they have rather distinct
graph structures, regardless of the common origin.

\begin{eqnarray*}
\begin{array}[b]{ c |  c c c c c c c c c }
\bullet & 1& 2& 3& 4& 5& 6& 7& 8& 9\\
\hline
1 & 9 & 7 & 3 & 9 & 9 & 7 & 3 & 7 & 3\\
2 & 9 & 7 & 3 & 9 & 9 & 7 & 3 & 7 & 3\\
3 & 9 & 7 & 3 & 9 & 9 & 7 & 3 & 7 & 3\\
4 & 8 & 6 & 2 & 8 & 8 & 6 & 2 & 6 & 2\\
5 & 5 & 4 & 1 & 5 & 5 & 4 & 1 & 4 & 1\\
6 & 8 & 6 & 2 & 8 & 8 & 6 & 2 & 6 & 2\\
7 & 8 & 6 & 2 & 8 & 8 & 6 & 2 & 6 & 2\\
8 & 5 & 4 & 1 & 5 & 5 & 4 & 1 & 4 & 1\\
9 & 5 & 4 & 1 & 5 & 5 & 4 & 1 & 4 & 1\\
\multicolumn{10}{c}{\mbox{ \tiny Natural central groupoid}}\\
\end{array}
& &
\epsfig{file=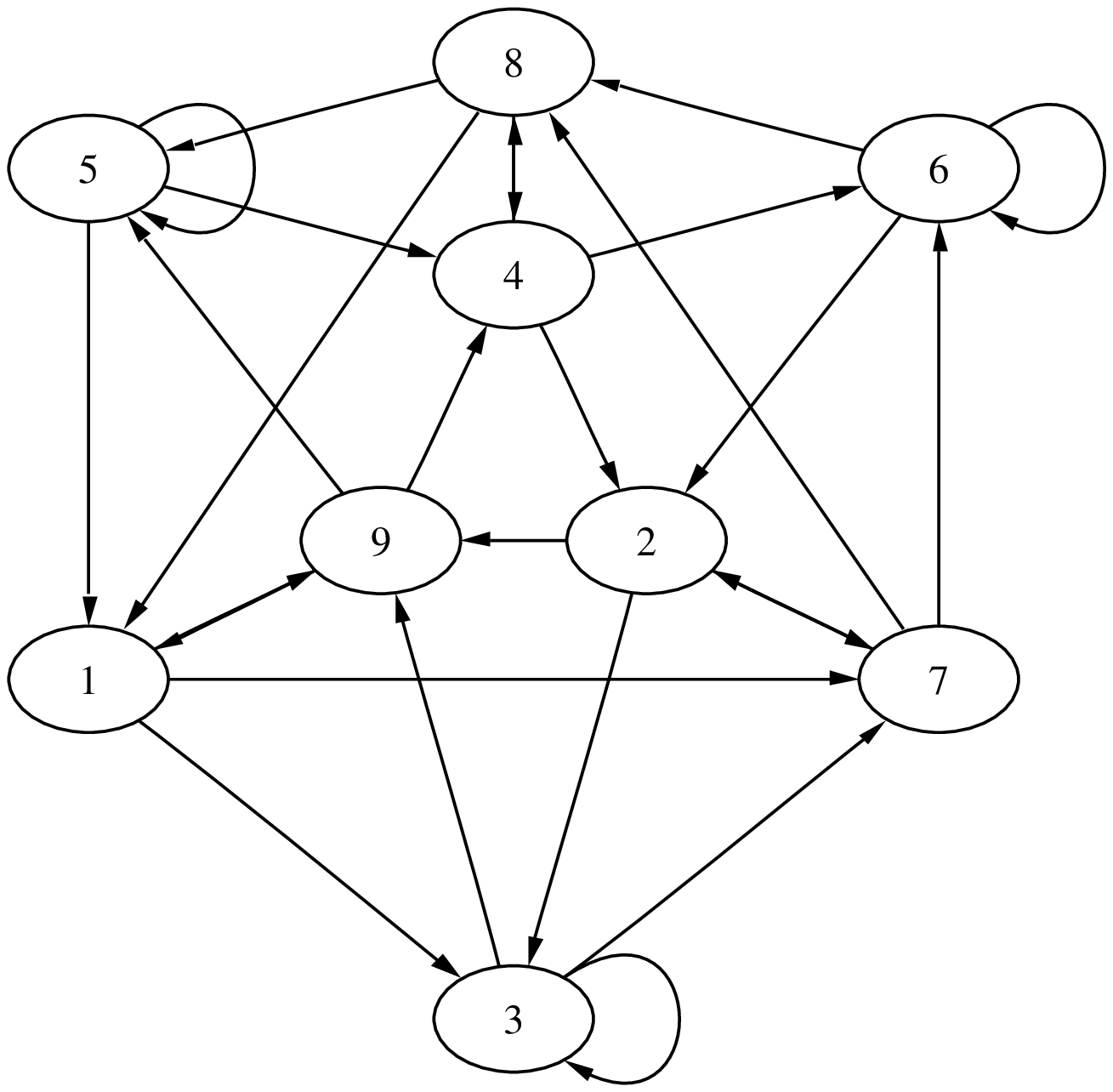,height=5.5cm}
\end{eqnarray*}
\begin{eqnarray*}
\begin{array}[b]{ c |  c c c c c c c c c  }
\bullet & 1& 2& 3& 4& 5& 6& 7& 8& 9\\
\hline
1 & 9 & 7 & 3 & 9 & 9 & 7 & 3 & 7 & 3\\
2 & 9 & 7 & 3 & 9 & 9 & 7 & 3 & 7 & 3\\
3 & 9 & 7 & 3 & 9 & 9 & 7 & 3 & 7 & 3\\
4 & 8 & 6 & 1 & 8 & 8 & 6 & 1 & 6 & 1\\
5 & 4 & 5 & 2 & 5 & 5 & 4 & 2 & 4 & 2\\
6 & 8 & 6 & 2 & 8 & 8 & 6 & 2 & 6 & 2\\
7 & 8 & 6 & 2 & 8 & 8 & 6 & 2 & 6 & 2\\
8 & 4 & 5 & 1 & 5 & 5 & 4 & 1 & 4 & 1\\
9 & 4 & 5 & 1 & 5 & 5 & 4 & 1 & 4 & 1\\
\multicolumn{10}{c}{\mbox{\tiny Number: 10 Lifting (1,9)(2,7)(4,8) }}\\
\end{array}& &
\epsfig{file=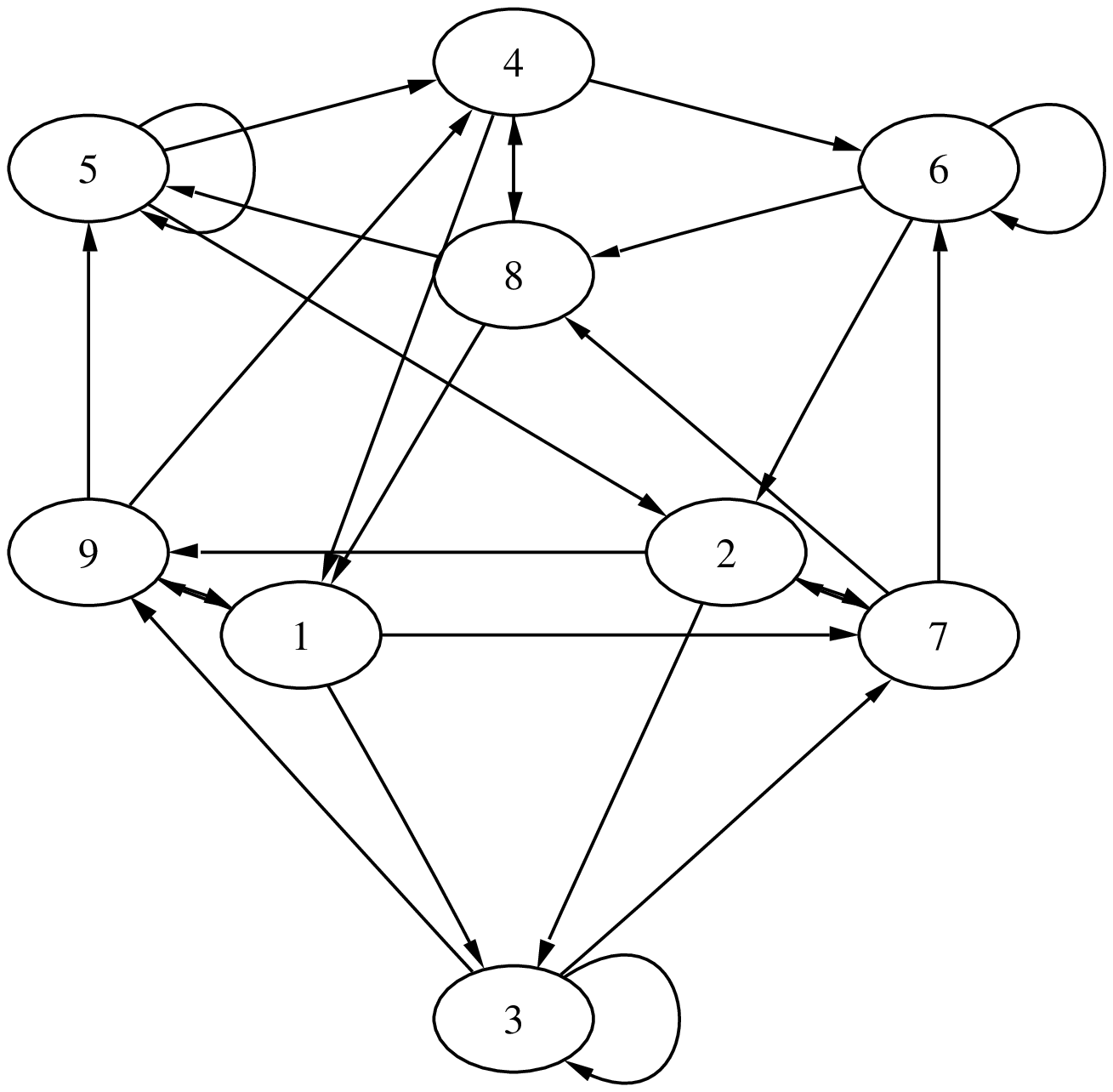,height=5.5cm}
\end{eqnarray*}

\begin{eqnarray*}
\begin{array}[b]{ c |  c c c c c c c c c  }
\bullet & 1& 2& 3& 4& 5& 6& 7& 8& 9\\
\hline
1 & 9 & 7 & 2 & 9 & 9 & 7 & 2 & 7 & 2\\
2 & 9 & 7 & 3 & 9 & 9 & 7 & 3 & 7 & 3\\
3 & 9 & 7 & 3 & 9 & 9 & 7 & 3 & 7 & 3\\
4 & 8 & 6 & 3 & 8 & 8 & 6 & 3 & 6 & 3\\
5 & 5 & 1 & 4 & 5 & 5 & 4 & 1 & 4 & 1\\
6 & 8 & 6 & 2 & 8 & 8 & 6 & 2 & 6 & 2\\
7 & 8 & 6 & 2 & 8 & 8 & 6 & 2 & 6 & 2\\
8 & 5 & 1 & 4 & 5 & 5 & 4 & 1 & 4 & 1\\
9 & 5 & 1 & 4 & 5 & 5 & 4 & 1 & 4 & 1\\
\multicolumn{10}{c}{\mbox{\tiny  Number: 36 Lifting (1,9)(2,7)(4,8) }}\\
\end{array}
& &
\epsfig{file=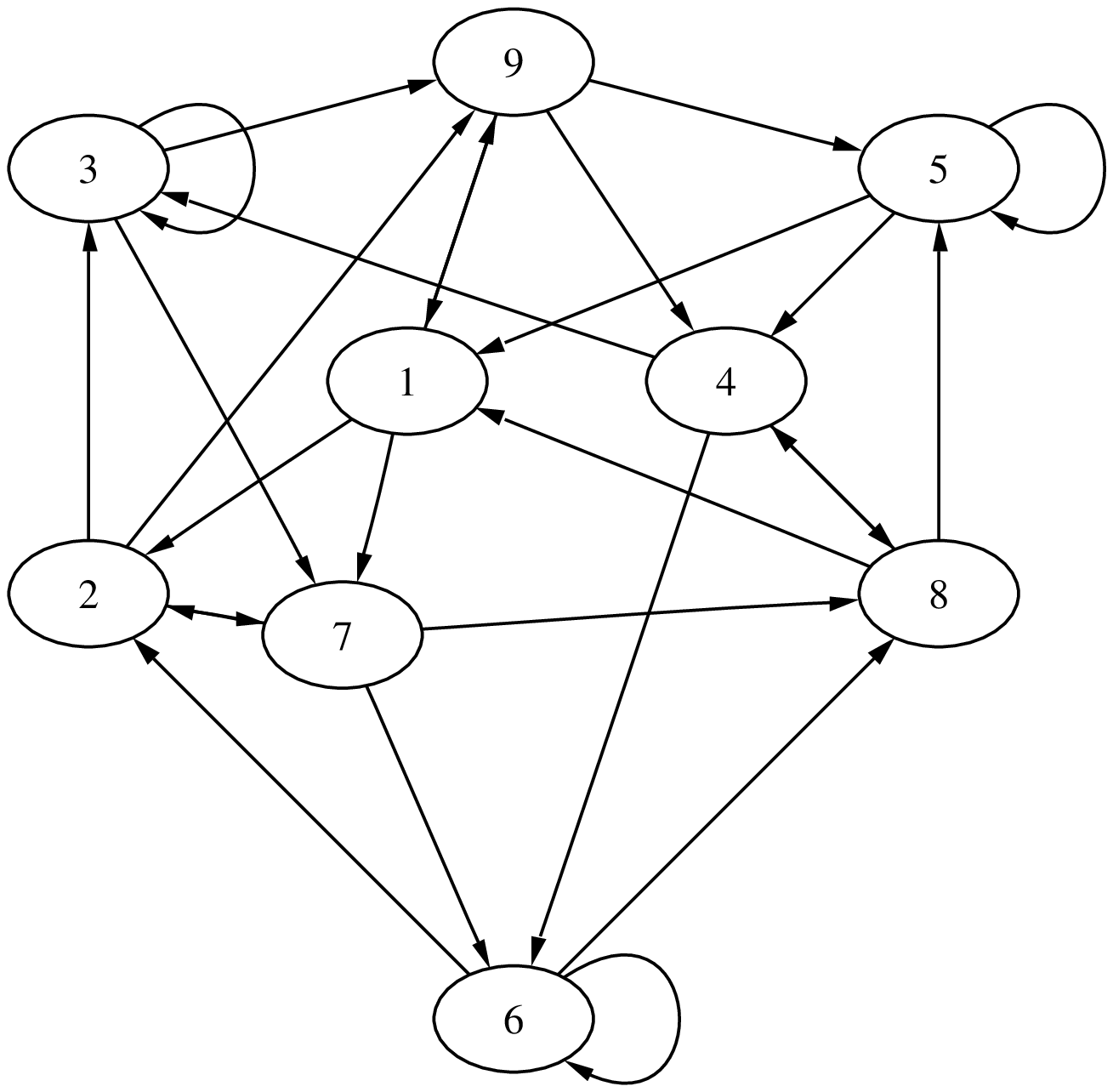,height=5.5cm}
\end{eqnarray*}
\begin{eqnarray*}
\begin{array}[b]{ c |  c c c c c c c c c}
\bullet & 1& 2& 3& 4& 5& 6& 7& 8& 9\\
\hline
1 & 9 & 7 & 3 & 9 & 9 & 7 & 3 & 7 & 3\\
2 & 9 & 4 & 3 & 9 & 9 & 4 & 3 & 4 & 3\\
3 & 9 & 7 & 3 & 9 & 9 & 7 & 3 & 7 & 3\\
4 & 8 & 6 & 2 & 2 & 8 & 6 & 8 & 6 & 2\\
5 & 5 & 4 & 1 & 5 & 5 & 4 & 1 & 4 & 1\\
6 & 8 & 6 & 2 & 2 & 8 & 6 & 8 & 6 & 2\\
7 & 8 & 6 & 2 & 2 & 8 & 6 & 8 & 6 & 2\\
8 & 5 & 7 & 1 & 5 & 5 & 7 & 1 & 7 & 1\\
9 & 5 & 4 & 1 & 5 & 5 & 4 & 1 & 4 & 1\\
\multicolumn{10}{c}{\mbox{\tiny Number: 105 Lifting(1,9)(2,4)(7,8) }}\\
\end{array}
& &
\epsfig{file=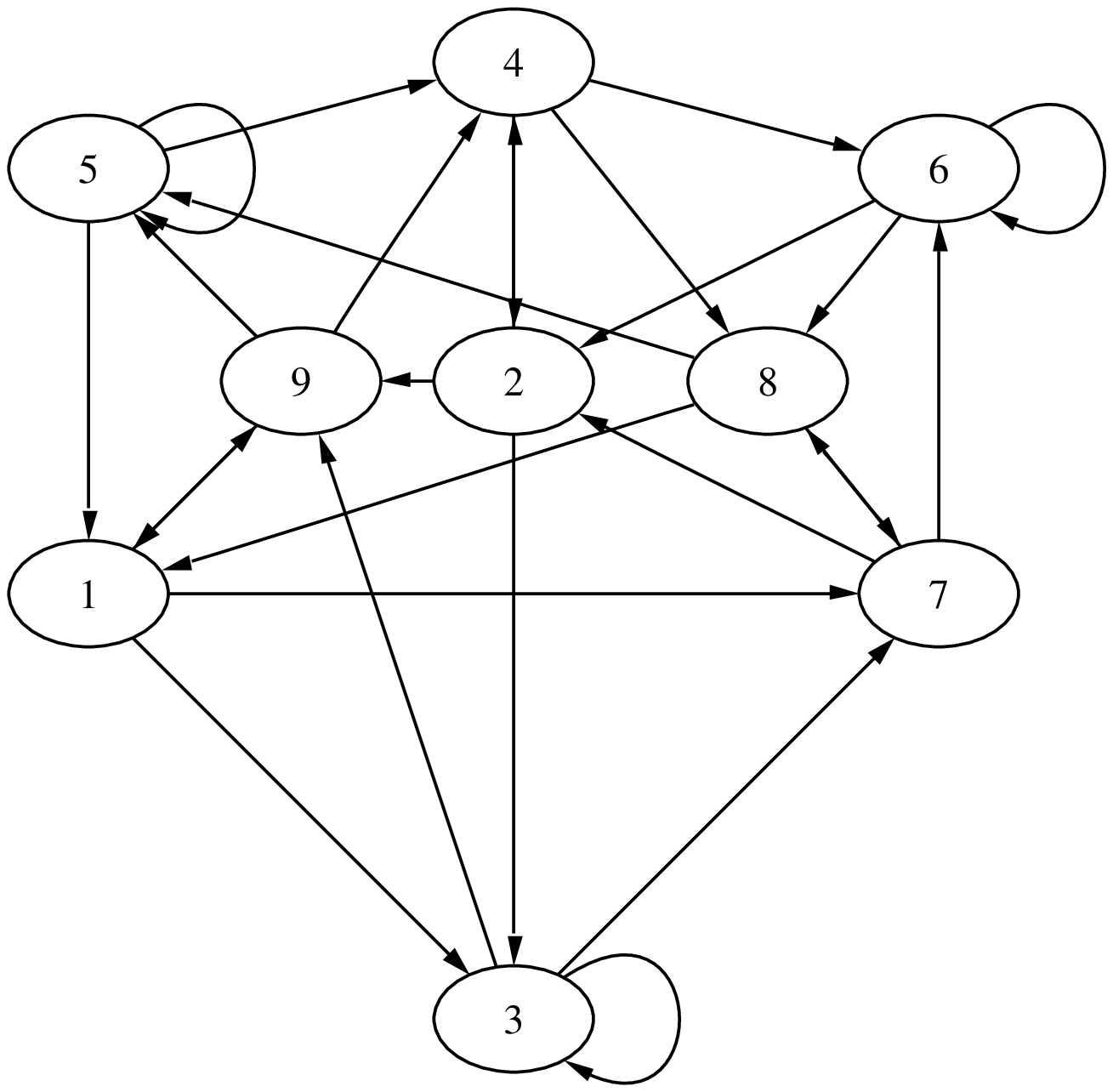,height=5.5cm}
\end{eqnarray*}
\begin{eqnarray*}
\begin{array}[b]{ c |  c c c c c c c c c  }
\bullet & 1& 2& 3& 4& 5& 6& 7& 8& 9\\
\hline
1 & 8 & 2 & 7 & 8 & 8 & 7 & 7 & 2 & 2\\
2 & 8 & 2 & 9 & 8 & 8 & 9 & 9 & 2 & 2\\
3 & 8 & 1 & 4 & 8 & 8 & 4 & 1 & 1 & 4\\
4 & 3 & 6 & 9 & 3 & 6 & 9 & 9 & 3 & 6\\
5 & 5 & 1 & 4 & 5 & 5 & 4 & 1 & 1 & 4\\
6 & 5 & 2 & 9 & 5 & 5 & 9 & 9 & 2 & 2\\
7 & 3 & 6 & 7 & 3 & 6 & 7 & 7 & 3 & 6\\
8 & 5 & 1 & 4 & 5 & 5 & 4 & 1 & 1 & 4\\
9 & 3 & 6 & 7 & 3 & 6 & 7 & 7 & 3 & 6\\
\multicolumn{10}{c}{\mbox{ \tiny Number: 118 Lifting (1,8)(3,4)(6,9) }}\\
\end{array}
& &
\epsfig{file=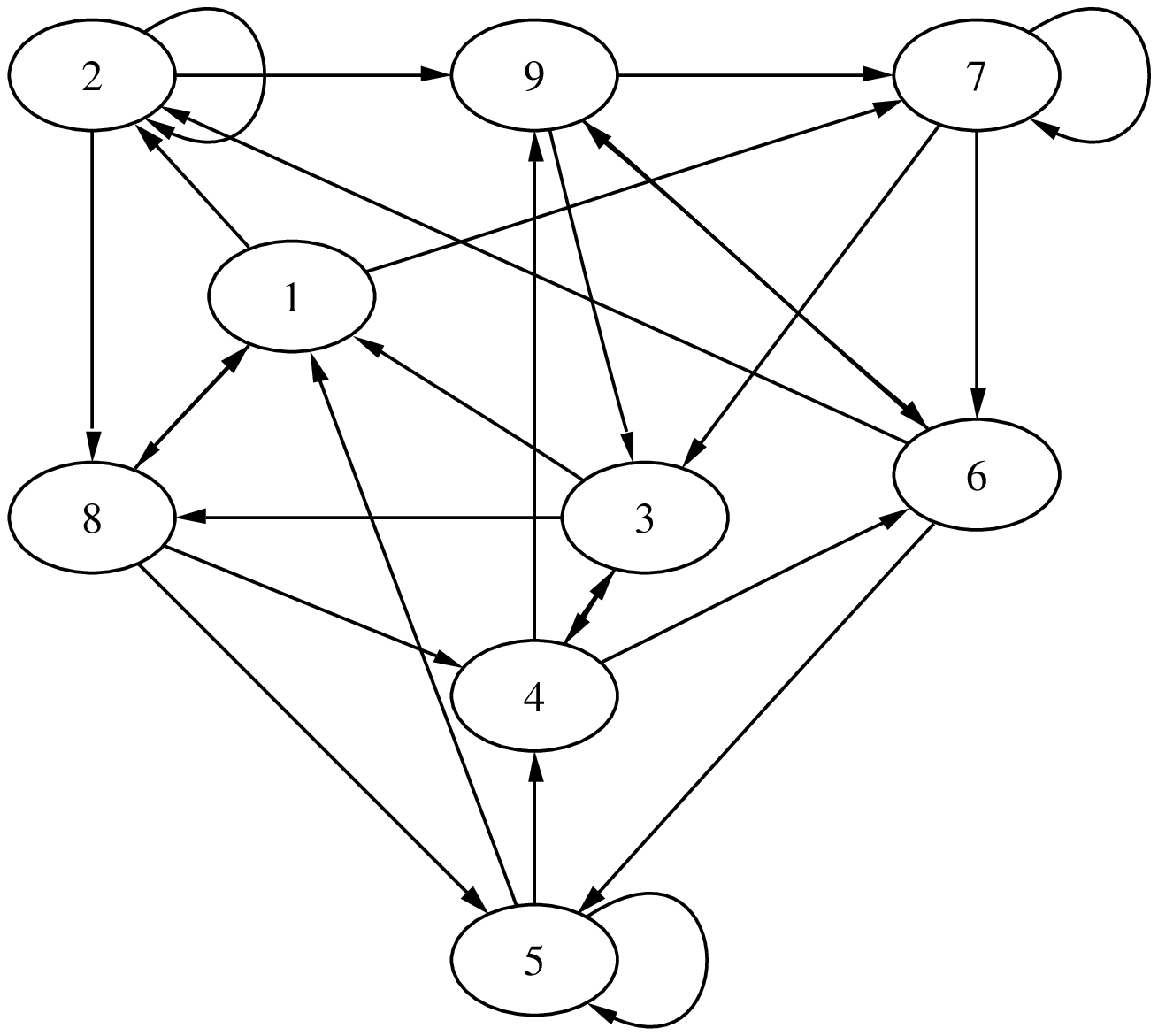,height=5.5cm}
\end{eqnarray*}
\begin{eqnarray*}
\begin{array}[b]{ c |  c c c c c c c c c  }
\bullet & 1& 2& 3& 4& 5& 6& 7& 8& 9\\
\hline
1 & 9 & 2 & 5 & 9 & 9 & 5 & 5 & 2 & 2\\
2 & 9 & 2 & 8 & 9 & 9 & 8 & 8 & 2 & 2\\
3 & 9 & 6 & 3 & 9 & 9 & 3 & 6 & 6 & 3\\
4 & 4 & 1 & 8 & 4 & 1 & 8 & 8 & 4 & 1\\
5 & 7 & 6 & 3 & 7 & 7 & 3 & 6 & 6 & 3\\
6 & 7 & 2 & 8 & 7 & 7 & 8 & 8 & 2 & 2\\
7 & 4 & 1 & 5 & 4 & 1 & 5 & 5 & 4 & 1\\
8 & 7 & 6 & 3 & 7 & 7 & 3 & 6 & 6 & 3\\
9 & 4 & 1 & 5 & 4 & 1 & 5 & 5 & 4 & 1\\
\multicolumn{10}{c}{\mbox{ \tiny Number: 118 Lifting (1,9)(5,7)(6,8)}}\\
\end{array}
& &
\epsfig{file=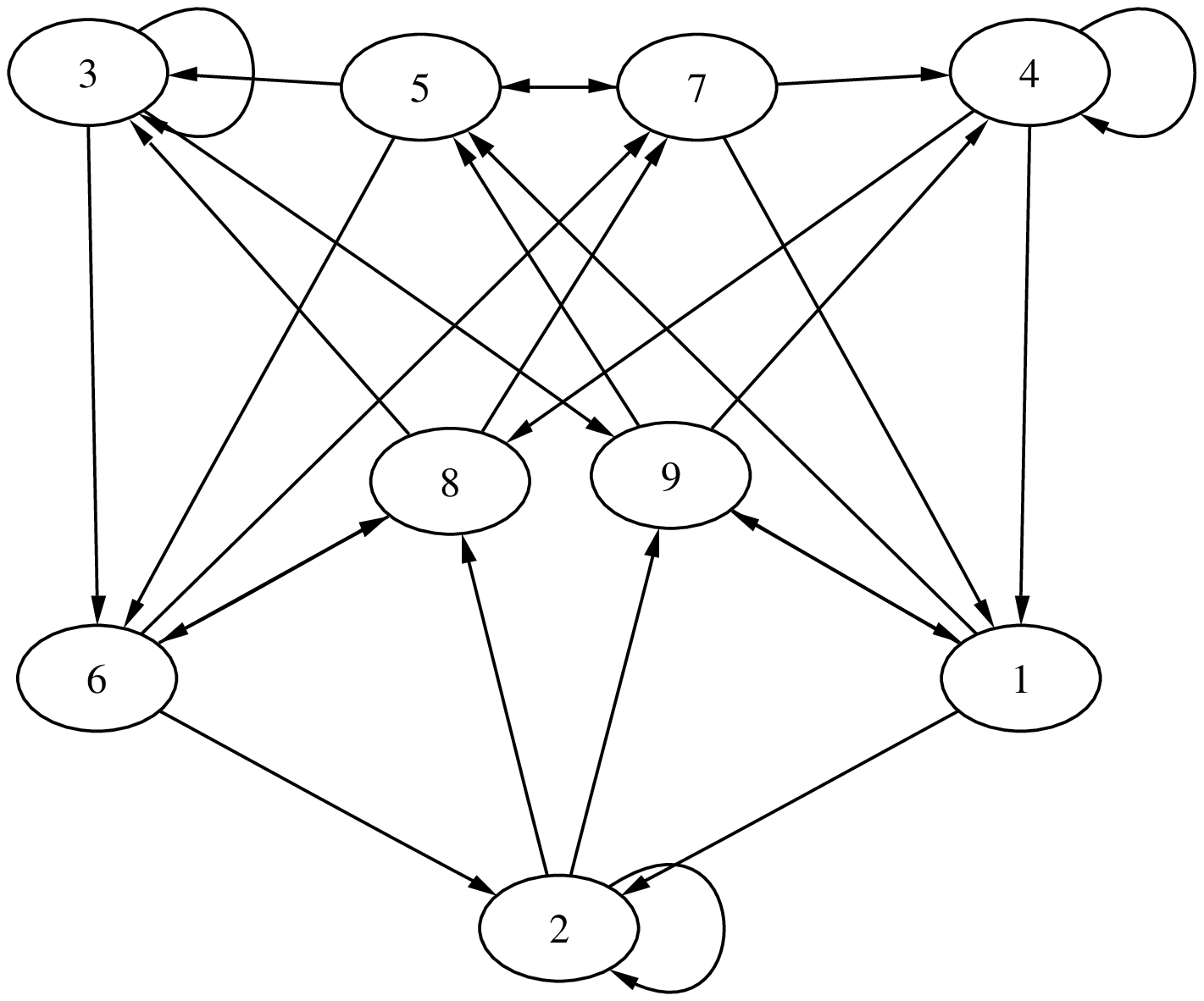,height=5.5cm}
\end{eqnarray*}

\section{Conclusion}

We have been able to use the techniques of orderly algorithms in order
to enumerate the central groupoids of order 4 and 9. It was apparently 
necessary to move to a more general, but combinatorially and computationally
more amenable structure, in order to be able to apply the
technique fruitfully. The combinatorial \rs\ was able to
be assembled in $n^2$ stages instead of $n^3$, which may have
suitably reduced the number of dead--end branches in the
generation tree.
The results connecting back to the
structure of central groupoids were then necessary in order
to be able to make the filtering process efficient.

The important results were the equivalence of idempotent \scbgs\ and
rectangular structures (Proposition \ref{prop_rect_to_scbg}) and the 
isomorphism of different liftings (Proposition \ref{propcon}).
Thus we were able to determine all central groupoids
and to remove isomorphs.

The use of combinatorial values in the orderly algorithm 
proved to be very important, reducing the search time from
over 360 minutes to less that 40.
It would be interesting to extend the results to an exhaustive list of 
central groupoids of order 16. 

The algorithms will need to be made more efficient in (probably)
several ways in order to make that search. It might be feasible
to define more combinatorial quantities on the graphs in order to lower 
the number of calls to {\tt nauty}. 
The technique that GRAPE uses to call {\tt nauty}, via text files
and remote invocation of the program, is possibly also not very 
efficient, as the
file system acts as a brake with every call.
The process of obtaining all possible extensions to a partial
\rs\ and then reducing modulo the automorphism group is possibly
acceleratable using some more efficient algorithm, perhaps
even an orderly algorithm.
Other analyses of the bottlenecks in this algorithm should
be made.

The implementations in GAP4, as well as the resulting collections
of examples, are available from the author electronically.

\section{Thanks}

The author would like to thank Leonard Soicher, who invited
him to speak about this work and gave the necessary impulse
to get the results written down. The author was supported in
part by project P15691 from the Austrian Federal FWF, the national
science finding body, as well as by several
ongoing grants from Stadt Linz, Land Ober\"osterreich and the
Austrian Federal BKA.Kunst.

\bibliographystyle{plain}
\bibliography{tims}

\end{document}